\documentclass[11pt]{amsart}

\usepackage{amsthm,amsmath,amssymb,amscd,amsfonts,latexsym}

\theoremstyle{plain}
\newtheorem{theorem}{Theorem}[section]

\newtheorem{fundamental fact}[theorem]{Fundamental Fact}

\newtheorem{corollary}[theorem]{Corollary}

\newtheorem{def-thm}[theorem]{Definition-Theorem}

\theoremstyle{definition}

\newtheorem{fact}[theorem]{Fact}

\newcommand{\PP}{\mathbb{P}}

\newcommand{\FF}{\mathbb{F}}

\newcommand{\OO}{{\mathcal O}}

\DeclareMathOperator{\kod}{kod}

\allowdisplaybreaks

\begin{document}

\title[Projectivized vector bundles and positive curvature]{On projectivized vector bundles and positive holomorphic sectional curvature}

\author{Angelynn Alvarez}
\address{Angelynn Alvarez. Department of Mathematics, University of Houston, 4800 Calhoun Road, Houston, TX 77204, USA} \email{{aalvarez@math.uh.edu}}

\author{Gordon Heier}
\address{Gordon Heier. Department of Mathematics, University of Houston, 4800 Calhoun Road, Houston, TX 77204, USA} \email{{heier@math.uh.edu}}

\author{Fangyang Zheng} \thanks{The third-named author is partially supported by a Simons Collaboration Grant}
\address{Fangyang Zheng. Department of Mathematics,
The Ohio State University, 231 West 18th Avenue, Columbus, OH 43210,
USA and Zhejiang Normal University, Jinhua, 321004, Zhejiang, China}
\email{{zheng.31@osu.edu}}

\subjclass[2010]{32L05, 32Q10, 32Q15, 53C55}
\keywords{Compact complex manifolds, K\"ahler metrics, positive holomorphic sectional curvature, positive scalar curvature, projectivized vector bundles}

\begin{abstract}
We generalize a construction of Hitchin to prove that, given any compact K\"ahler manifold $M$ with positive holomorphic sectional curvature and any holomorphic vector bundle $E$ over $M$, the projectivized vector bundle ${\mathbb P}(E)$ admits a K\"ahler metric with positive holomorphic sectional curvature.
\end{abstract}

\maketitle

\section{Introduction}

Let us denote by ${\mathcal H}_n$ (resp.\ ${\mathcal H}_n'$) the set of all $n$-dimensional projective manifolds (resp.\ compact complex manifolds) which admit K\"ahler metrics with positive holomorphic sectional curvature $H$. It is a long-standing goal in complex geometry to understand these sets, which are conjectured to agree. As far as we know, this agreement is still open for $n\geq 3$. \par

First, let us recall some known facts. For any $M$ in ${\mathcal H}_n'$, a theorem of Tsukamoto \cite{Tsukamoto_1957} states that $M$ is simply-connected. Furthermore, by a theorem of Berger \cite{Berger}, the scalar curvature of any K\"ahler metric on $M$ is the average of the holomorphic sectional curvature. Therefore, $H>0$ implies that the scalar curvature is positive. It was proven by Kobayashi and Wu \cite[Corollary 2]{Kobayashi_Wu} that consequently the Kodaira dimension $\kod(M) = -\infty $. In recent work by the second named author and Wong \cite{Heier_Wong_1509} it has been proven that any $M$ in ${\mathcal H}_n$ is rationally connected, i.e., given any two points in $M$, there exist a connected chain of rational curves containing both points. \par

On the existence side, it is well-known that the Fubini-Study metric on ${\mathbb P}^2$ has constant positive $H$, so it is in ${\mathcal H}_2$. In \cite{Hitchin}, Hitchin showed that any Hirzebruch surface $\FF_a=\PP(\OO_{\PP^1}(a) \oplus  \OO_{\PP^1})$, $a\in \{0,1,2,\ldots\}$, is also in ${\mathcal H}_2$. For all other rational surfaces, including the del Pezzo surfaces, it is not known to the best of our knowledge whether or not they admit any K\"ahler metric with positive $H$. In general, there are only few known elements of ${\mathcal H}_n'$. Products of manifolds with positive $H$ constitute a further obvious class of elements. Also, on K\"ahler $C$-spaces (i.e., simply-connected, compact, homogeneous K\"ahler manifolds), the standard K\"ahler-Einstein metric has positive holomorphic sectional curvature, so they all belong to ${\mathcal H}_n$.\par

The main purpose of this note is to generalize Hitchin's construction on Hirzebruch surfaces to arbitrary projectivized vector bundles over a compact K\"ahler manifold with positive holomorphic sectional curvature and thus provide a new class of examples of manifolds in ${\mathcal H}_n'$ as stated in the following theorem.
\begin{theorem} \label{main_thm} Let $M$ be a compact K\"ahler manifold with positive holomorphic sectional curvature. Let $E$ be a holomorphic vector bundle over $M$ and $P={\mathbb P}(E)$ the projectivization of $E$. Then $P$ admits a K\"ahler metric with positive holomorphic sectional curvature.
\end{theorem}
If we replace ${\mathbb P}(E)$ by the Grassmannian bundle $G_k(E)$ of all $k$-dimensional subspaces of the fibers of $E$ (here $k$ is any positive integer less than the rank of $E$), then it is clear from the proof of the theorem that $G_k(E)$ also has positive holomorphic sectional curvature, and it can be holomorphically and isometrically embedded into ${\mathbb P}(\Lambda ^kE)$.  \par

Recall that in \cite{Hitchin}, the proof of the existence of a metric of positive holomorphic sectional curvature on Hirzebruch surfaces is a stepping stone in the proof of the existence of metrics of positive scalar curvature on generic rational surfaces. As part of the argument, the following fundamental fact is established there as \cite[Corollary (5.18)]{Hitchin}.

\begin{fact}
Let $X$ be a compact K\"ahler manifold of dimension $>2$ with everywhere positive scalar curvature. Suppose we blow up any point $p\in X$ to obtain $\hat X$; then $\hat X$ admits a K\"ahler metric of positive scalar curvature.
\end{fact}

Based on this fact and Berger's theorem, we observe the following immediate corollary to Theorem \ref{main_thm}.

\begin{corollary}
Let $P$ be as in Theorem \ref{main_thm}. Assume that the dimension of $P$ is greater than $2$. Let $Y$ be obtained from $P$ by a finite sequence of blow-ups of points. Then $Y$ carries a K\"ahler metric of positive scalar curvature.
\end{corollary}

In general, it is reasonable to conjecture that for any holomorphic fiber bundle with total space $X$, if both the fiber and the base admit K\"ahler metrics with positive holomorphic sectional curvature, then so does $X$. However, at this point we do not know how to construct a natural metric on $X$ from the given metrics on the fiber and the base in a way that enables us to compute the holomorphic sectional curvature. In the case of negative holomorphic sectional curvature, an even more general result of this nature was obtained by Cheung \cite[Theorem 1]{Cheung}.

We conclude this introduction by remarking that in the paper \cite{ACH} an explicit analysis of the pinching constants for Hitchin's metrics on Hirzebruch surfaces was conducted. We leave such an explicit analysis for the present higher-dimensional case to a future occasion.

\section{Proof of Theorem \ref{main_thm}}

In this section, we will prove Theorem \ref{main_thm}. Let $(M,g)$ be a compact $n$-dimensional K\"ahler manifold with positive holomorphic sectional curvature. Let $(E,h)$ be a holomorphic vector of rank $r+1$ equipped with a Hermitian metric. Denote by $P={\mathbb P}(E)$ the projectivization of $E$, namely, for any $x\in M$, the fiber $P_x$ is simply the projective space ${\mathbb P}(E_x)$ of the fiber $E_x$. In other words, $P_x$ consists of the equivalence classes $[w]$ where $w$ is any non-zero vector in $E_x$. Note that this is the differential-geometric notation -- the algebro-geometric notation of ${\mathbb P}(E)$ would be the ${\mathbb P}(E^{\ast })$ here. \par

As essentially observed in \cite[(4.1)]{Hitchin}, the metrics $g$ and $h$ naturally induce a closed $(1,1)$-form on $P$:
$$ \omega_G = \lambda \ \pi^{\ast }(\omega_g) + \sqrt{-1}\  \partial \overline{\partial } \log h(v,\overline{v}), $$
where $\omega_g$ is the K\"ahler form of $g$, $\pi : P\rightarrow M$ the projection map, and $(x,[v])$ is a moving point in $P$. Since the restriction of the second term of the right hand side on a fiber of $\pi$ is just the Fubini-Study metric, we know that for $\lambda$ sufficiently large, $\omega_G$ is positive definite everywhere so $G=G_{\lambda}$ becomes a K\"ahler metric on $P$. \par

We claim that there is some constant $\lambda_0>0$ which depends on $g$ and $h$, such that for any $\lambda \geq \lambda_0$, the metric $G=G_{\lambda}$ has positive holomorphic sectional curvature. This will complete the proof of our theorem.\par

To prove this, we need to compute the curvature of the metric $G$. Fix any $p=(x_0,[w])$ in $P$. Without loss of generality, we may assume that $|w|=1$. Let $(z_1, \ldots , z_n)$ be a local holomorphic coordinate centered at $x_0$ which is normal with respect to $g$. That is, $x_0=(0,\ldots , 0)$, and under the coordinates,
$$g_{i\overline{j}}(0)=\delta_{i\overline{j}}, \ \ \ dg_{i\overline{j}}(0)=0$$
 for any $1\leq i,j\leq n$. Write $\Theta^h$ for the curvature of $(E,h)$. Then by a constant unitary change of $z$ if necessary, we may assume that the $(1,1)$-form $\Theta^h_{w\overline{w}}$ at $x_0$ is diagonal, namely,
$$ \Theta^h_{w\overline{w}} = \sum_{i=1}^n \xi_i \ dz_i \wedge d\overline{z}_i.$$
Next let us choose a local holomorphic frame $\{ e_0, e_1, \ldots , e_r\}$ of $E$ near $x_0$. Write $h_{\alpha \overline{\beta}}$ for $h(e_{\alpha }, \overline{e_{\beta }})$. We may assume that at the origin, $e_0(0)=w$, and
$$ h_{\alpha \overline{\beta}} (0)=\delta_{\alpha \beta}, \ \ dh_{\alpha \overline{\beta}} (0)=0.$$
It is also easy to see that we may further assume that  $\partial_i\partial_k h_{\alpha \overline{\beta}}(0)=0$. Let us write
$$ v= e_0(z) + \sum_{\alpha =1}^r t_{\alpha } e_{\alpha}(z) .$$
Then $(z,t)$ becomes a local holomorphic coordinate in $P$ centered at $p=(0,0)$. We will use the index convention of $z_i$, $t_{\alpha}$, and with index after comma denoting the partial derivatives. Because of the sesquilinearity of $h$, we may naturally write $h_{\alpha \overline{v}}$ instead of $h_{v\overline{v},\alpha}$ and  $h_{v\overline{\beta}}$ instead of $h_{v\overline{v},\overline{\beta}}$. We have
\begin{eqnarray*}
G_{i\overline{j}}  & = & \lambda \ g_{i\overline{j}} +  \frac{1} { h_{v\overline{v}} } h_{v\overline{v}, i\overline{j}} - \frac{1}{ (h_{v\overline{v}})^2 } h_{v\overline{v}, i}  \ h_{v\overline{v}, \overline{j} } \\
G_{i\overline{\beta} }  & = & \frac{1}{h_{v\overline{v}}}  h_{v\overline{\beta}, i} -  \frac{1}{ (h_{v\overline{v}})^2 }  h_{v\overline{v}, i} \ h_{v\overline{\beta}} \\
 G_{\alpha\overline{j}}  & = & \frac{1}{h_{v\overline{v}}}  h_{\alpha \overline{v}, \overline{j}} -  \frac{1}{ (h_{v\overline{v}})^2 } h_{v\overline{v}, \overline{j}} \  h_{\alpha \overline{v} } \\
 G_{\alpha \overline{\beta} }  & = & \frac{1}{h_{v\overline{v}}}  h_{\alpha \overline{\beta}, i} -  \frac{1}{ (h_{v\overline{v}})^2 } h_{\alpha \overline{v}} \ h_{v\overline{\beta}}.
 \end{eqnarray*}
So at the origin $p$, we have $G_{i\overline{j}} (0) = (\lambda -\xi_i) \ \delta_{ij} $, \ $G_{i\overline{\beta}} (0) = G_{\alpha \overline{j}}(0)=0$, and $G_{\alpha \overline{\beta}}(0)=\delta_{\alpha\beta}$, for any $1\leq  i,j\leq n$ and any $1\leq \alpha , \beta \leq r$. Next, we compute the first derivatives of $G$:
\begin{eqnarray*}
G_{i\overline{j} , k}  & = & \lambda \ g_{i\overline{j},k } +  \frac{1} { h_{v\overline{v}} } h_{v\overline{v}, i\overline{j}k} - \frac{1}{ (h_{v\overline{v}})^2 } h_{v\overline{v}, k}  \ h_{v\overline{v}, i\overline{j} } -  \frac{1}{ (h_{v\overline{v}})^2 } h_{v\overline{v}, ik}  \ h_{v\overline{v}, \overline{j} } \\
&  & -  \frac{1}{ (h_{v\overline{v}})^2 } h_{v\overline{v},i}  \ h_{v\overline{v}, \overline{j}k } +  \frac{2}{ (h_{v\overline{v}})^3 } h_{v\overline{v}, k}  \ h_{v\overline{v}, i }\  h_{v\overline{v}, \overline{j} } \\
G_{i\overline{\beta} ,k }  & = & \frac{1}{h_{v\overline{v}}}  h_{v\overline{\beta}, ik} -  \frac{1}{ (h_{v\overline{v}})^2 }  h_{v\overline{v}, k} \ h_{v\overline{\beta},i}   -  \frac{1}{ (h_{v\overline{v}})^2 }  h_{v\overline{v}, i} \ h_{v\overline{\beta},k}  \\
& &  -  \frac{1}{ (h_{v\overline{v}})^2 }  h_{v\overline{v}, ik} \ h_{v\overline{\beta}}   +   \frac{2}{ (h_{v\overline{v}})^3 }  h_{v\overline{v}, i} \  h_{v\overline{v}, k} \ h_{v\overline{\beta}} \\
 G_{\alpha\overline{j},k}  & = & \frac{1}{h_{v\overline{v}}}  h_{\alpha \overline{v}, \overline{j}k} -   \frac{1}{ (h_{v\overline{v}})^2 } h_{v\overline{v}, k} \  h_{\alpha \overline{v}, \overline{j} }   -   \frac{1}{ (h_{v\overline{v}})^2 } h_{v\overline{v}, \overline{j}} \  h_{\alpha \overline{v},k } \\
 & & -  \frac{1}{ (h_{v\overline{v}})^2 } h_{v\overline{v}, \overline{j}k } \  h_{\alpha \overline{v} }  + \frac{2}{ (h_{v\overline{v}})^3 }   h_{\alpha \overline{v} } h_{v\overline{v}, \overline{j} } \  h_{v\overline{v},k }  \\
 G_{\alpha \overline{\beta} ,k}  & = & \frac{1}{h_{v\overline{v}}}  h_{\alpha \overline{\beta}, k} -  \frac{1}{ (h_{v\overline{v}})^2 } h_{\alpha \overline{\beta}} \ h_{v\overline{v},k}  - \frac{1}{ (h_{v\overline{v}})^2 } h_{\alpha \overline{v},k } \ h_{v\overline{\beta}}\\
 & & - \frac{1}{ (h_{v\overline{v}})^2 } h_{\alpha \overline{v}} \ h_{v\overline{\beta},k} + \frac{2}{ (h_{v\overline{v}})^3 } h_{\alpha \overline{v}} \ h_{v\overline{\beta}} \  h_{v\overline{v},k} \\
  G_{\alpha\overline{j},\gamma }  & = &  -   \frac{1}{ (h_{v\overline{v}})^2 } h_{\gamma \overline{v}} \  h_{\alpha \overline{v}, \overline{j} }   -   \frac{1}{ (h_{v\overline{v}})^2 } h_{\gamma \overline{v}, \overline{j}} \  h_{\alpha \overline{v} }  + \frac{2}{ (h_{v\overline{v}})^3 }   h_{\alpha \overline{v} } \  h_{\gamma \overline{v} } \ h_{v\overline{v}, \overline{j} }   \\
 G_{\alpha \overline{\beta} ,\gamma }  & = &  -  \frac{1}{ (h_{v\overline{v}})^2 } h_{\alpha \overline{\beta}} \ h_{\gamma \overline{v}}  - \frac{1}{ (h_{v\overline{v}})^2 } h_{\alpha \overline{v} } \ h_{\gamma \overline{\beta}} + \frac{2}{ (h_{v\overline{v}})^3 } h_{\alpha \overline{v}} \ h_{\gamma \overline{v}} \ h_{v\overline{\beta}}.
 \end{eqnarray*}
We also have $G_{i\overline{j},\alpha } = G_{\alpha \overline{j}, i}$ and  $G_{i\overline{\beta},\alpha } = G_{\alpha \overline{\beta}, i}$ by the K\"ahlerness of $G$. At the origin $p$, we have $h_{v\overline{v}}(0)=|w|^2 =1$, $h_{v\overline{\beta}}(0)=h_{\alpha \overline{v}}(0)=0$, and all first order derivatives of $h$ are zero, so by taking another derivative and evaluate at $0$, we get the following at the point $p$:
\begin{eqnarray*}
G_{i\overline{j},k\overline{l}} & = & \lambda g_{i\overline{j},k\overline{l}} + h_{v\overline{v},i\overline{j}k\overline{l}} - h_{v\overline{v},i\overline{j}} h_{v\overline{v},k \overline{l}} - h_{v\overline{v},i\overline{l}} h_{v\overline{v},k\overline{j}} \\
G_{i\overline{j},k\overline{\beta}} & = & h_{v\overline{\beta},i\overline{j}k} \\
G_{i\overline{j},\alpha \overline{\beta}} & = & h_{\alpha \overline{\beta},i\overline{j}} - h_{\alpha\overline{\beta}}  h_{v\overline{v},i\overline{j}} \\
G_{\alpha\overline{j},\gamma \overline{l}} & = & 0 \\
G_{\alpha\overline{\beta},\gamma \overline{j}} & = & 0\\
G_{\alpha\overline{\beta},\gamma \overline{\delta}} & = & - h_{\alpha\overline{\beta}} h_{\gamma \overline{\delta}} - h_{\alpha\overline{\delta}} h_{\gamma\overline{\beta}}.
\end{eqnarray*}
Now we are ready to compute the holomorphic sectional curvature of $G$ at $p$. Let $0\neq V=X+U$ be a type $(1,0)$ tangent vector at $p\in P$, where
$X= \sum_{i=1}^n x_i \frac{\partial }{\partial z_i}$ and $U = \sum_{\alpha =1}^r u_{\alpha} \frac{\partial }{\partial t_{\alpha }} $.  Denote by $R$, $R^g$, and $R^h$  the curvature tensors of $G$, $g$, and $h$, respectively.  We have
 \begin{eqnarray*}
 R_{V \overline{V}V \overline{V}} & = & R_{X \overline{X}X \overline{X}} + 4R_{X \overline{X}U \overline{U}} + R_{U \overline{U}U\overline{U}} \\
  & & + \ 2Re\{ R_{X \overline{U}X \overline{U}}+2R_{X \overline{X}X \overline{U}}+2R_{U \overline{U}U \overline{X}} \}
  \end{eqnarray*}
by the symmetry of the curvature tensor. At $p$ the matrix of $G$ is diagonal, so we have
 \begin{eqnarray*}
 R_{V \overline{V}V\overline{V}} & = & - G_{V\overline{V},V\overline{V}} + \sum_{a=1}^{n+r} \frac{1}{G_{a\overline{a}} } |G_{V\overline{a},V}|^2 \\
 & \geq & - G_{V\overline{V},V\overline{V}} \\
 & = & -G_{X \overline{X},X \overline{X}} -4G_{X \overline{X},U \overline{U}} -G_{U \overline{U},U\overline{U}} \\
 & & - \ 2Re\{ G_{X \overline{U},X \overline{U}}+2G_{X \overline{X},X \overline{U}}+2G_{U \overline{U},U \overline{X}} \} \\
 & = &  -G_{X \overline{X},X \overline{X}} -4G_{X \overline{X},U \overline{U}} -G_{U \overline{U},U\overline{U}}  - 4Re\{G_{X \overline{X},X \overline{U}} \} \\
 & = & (\lambda R^g_{X \overline{X}X \overline{X}} - h_{v\overline{v},X \overline{X}X \overline{X}} + 2 (R^h_{v\overline{v}X \overline{X}} )^2 )  +4(R^h_{U\overline{U}X\overline{X}} \\
 & &   -\ |U|^2 R^h_{v\overline{v}X\overline{X}}) + 2|U|^4 - 4Re\{ h_{v\overline{U}, X\overline{X}X} \}.
 \end{eqnarray*}
Let $H_0>0$ be the minimum of holomorphic sectional curvature of $g$ over $M$, and choose a constant $C>0$ such that $|R^h|$ and $|\nabla R^h|$ are both bounded by $C$, then the above computation leads to
$$ R_{V \overline{V}V\overline{V}}  \geq  (\lambda H_0 - C) |X|^4 - 8C |X|^2 |U|^2 + 2|U|^4 - 4C |X|^3 |U|.$$
Clearly, if $\lambda $ is sufficiently large, this quantity will be positive when $X$ and $U$ are not both zero. This completes the proof of Theorem \ref{main_thm}.

\end{document}